\documentclass[10pt,leqno]{amsart}

\usepackage {yfonts}
\usepackage [T1] {fontenc}

\usepackage{xcolor}
\usepackage{textcomp}

\usepackage{enumitem}
\usepackage{soul}

\numberwithin{equation}{section}
\newtheorem{theorem}{Theorem}[section]
\newtheorem{prop}[theorem]{Proposition}
\newtheorem{lem}[theorem]{Lemma}

\newtheorem{rem}[theorem]{Remark}

\newcommand{\R}{\mathbb{R}}

\newcommand{\N}{\mathbb{N}}

\newcommand{\F}{\mathcal{F}}
\theoremstyle{plain}

\def\pn{\hfil\par\noindent}
\def\be{\begin{enumerate}}
\def\ds{\displaystyle }
\def\ee{\end{enumerate}}

\def\beqq{\begin{eqnarray*}}
\def\eeqq{\end{eqnarray*}}

\def\buildo#1\over#2{\mathrel{\mathop{\null#2}\limits^{#1}}}
\def\buildu#1\under#2{\mathrel{\mathop{\null#2}\limits_{#1}}}

\title{Long Time Decay of Leray Solution of 3D-NSE With Damping}
\author{Mongi Blel}
   \address{King Saud University, College of Sciences, Department of Mathematics,   Kingdom of Saudi Arabia}
   \email{mblel@ksu.edu.sa, jamelbenameur@gmail.com}
\author{Jamel Benameur}

\date{\today}


\begin{document}

\maketitle

\begin{abstract}
In \cite{CJ}, the authors show that the Cauchy problem of the Navier-Stokes equations with damping $\alpha|u|^{\beta-1}u(\alpha>0,\;\beta\geq1)$ has global
weak solutions in $L^2(\R^3)$. In this paper, we prove the uniqueness, the continuity in $L^2$ for $\beta>3$, also the large
time decay is proved for $\beta\geq\frac{10}3$. Fourier analysis and standard
techniques are used.
\end{abstract}

\section{\bf Introduction}\ \\
In this paper we study the global existence of weak solution to the modified incompressible Navier-Stokes equations in $\R^3$

\begin{equation}\label{$S$}
 \left\{ \begin{matrix}
     \partial_t u
 -\nu\Delta u+ u.\nabla u  +\alpha|u|^{\beta-1}u = -\nabla p &\hbox{ in } \mathbb R^+\times \mathbb R^3\\
     {\rm div}\, u = 0 \hfill&\hbox{ in } \mathbb R^+\times \mathbb R^3\\
    u(0,x) =u^0(x) \hfill&\hbox{ in }\mathbb R^3,\\
    \alpha>0,\;\beta>1\hfill&.
  \end{matrix}\right. \tag{$NSD$}
\end{equation}
 where $u=u(t,x)=(u_1,u_2,u_3)$, $p=p(t,x)$ denote respectively the unknown velocity and the pressure of the fluid at the point $(t,x)\in \mathbb R^+\times \mathbb R^3$, $\nu $ is the viscosity of fluid  and $u^0=(u_1^0(x),u_2^0(x),u_3^0(x))$ the initial given velocity.
 The damping is from the resistance to the motion of the flow. It describes various physical situations such as porous media flow, drag or friction effects, and some dissipative mechanisms (see \cite{BD,BDC,H,HP} and references therein).
  The fact that ${\rm div}\,u = 0$, allows to write the term $(u.\nabla u):=u_1\partial_1 u+u_2\partial_2 u+u_3\partial_3u$ in the following form
$ {\rm div}\,(u\otimes u):=({\rm div}\,(u_1u),{\rm div}\,(u_2u),{\rm div}\,(u_3u)).$
 If the initial   velocity $u^0$ is quite regular, the divergence free condition determines the pressure $p$.\\
  In order to simplify the calculations and the proofs of our results, we consider the viscosity  unitary (i.e. $\nu=1$).

  The global existence of weak solution of initial value problem of the classical incompressible Navier-Stokes were proved by Leray and Hopf (see \cite{Hopf}-\cite{Leray}) long before. The uniqueness remains an open problem for the dimensions $d\geq3$.\\
  The polynomial damping $\alpha|u|^{\beta-1}u$ is studied in \cite{CJ} by Cai and Jiu, where they proved the global  existence of weak solution in
 $$L^\infty(\R^+,L^2(\R^3))\cap L^2(\R^+,\dot H^1(\R^3))\cap L^{\beta+1}(\R^+,L^{\beta+1}(\R^3)).$$

 The purpose of this paper is to study the uniqueness, continuity and large time decay of the global solution of the incompressible Navier-Stokes equations with damping $(NSD)$. We recall that in \cite{CJ} employ the Galerkin approximations to construct the global solution of $(NSD)$ with $\beta\geq1$. But, in our case we use Friedrich method to prove the continuity and the uniqueness of such solution for $\beta>3$. The study of large time decay is studied for $\beta\geq \frac{10}3$. Precisely, our main result is the following:
\begin{theorem}\label{th2}\pn
 Let $\beta>3$ and $u^0\in L^2(\mathbb R^3)$ be a divergence free vector fields, then there is a unique global solution of $(NSD)$:
$u\in C_b(\R^+,L^2(\mathbb R^3)\cap L^2(\R^+,\dot H^1(\mathbb R^3))\cap L^{\beta+1}(\R^+,L^{\beta+1}(\R^3))$. Moreover, for all $t\geq0$
\begin{equation}\label{eqth2-1}\|u(t)\|_{L^2}^2+2\int_0^t\|\nabla u(s)\|_{L^2}^2ds +2\alpha\int_0^t\|u(s)|_{L^{\beta+1}}^{\beta+1}ds
\leq \|u^0\|_{L^2}^2.
\end{equation}
Moreover, if $\beta\geq\frac{10}3$ we have
\begin{equation}\label{eqth2-2}
\limsup_{t\to \infty}\|u(t)\|_{L^2}=0.
\end{equation}
\end{theorem}
\begin{rem}\pn
  In theorem \ref{th2}, the inequality (\ref{eqth2-1}) is proved in
  \cite{CJ}. The new parts of this theorem is the uniqueness, the continuity of the
  global solution in $L^2(\R^3)$ and the asymptotic result (\ref{eqth2-2}).
\end{rem}

\section{\bf Notations and Preliminary Results}
For a function  $f\colon\R^3\to\bar\R$    and $R>0$, the
Friedritch operator $J_R$  is defined by:
  $\ds J_R(D)f=\F^{-1}(\chi_{B_R} \widehat{f}),$
where $B_R$,   the ball of center $0$ and radius $R$. If $L^2_\sigma(\R^3)$ is the space of divergence-free vector fields in $L^2 (\R^3)$, the
 Leray projector $\mathbb P\colon (L^2(\R^3))^3\to (L^2_\sigma(\R^3))^3$ is  defined by:
$$\mathcal F(\mathbb P f)=\widehat{f}(\xi)-(\widehat{f}(\xi).\frac{\xi}{|\xi|})\frac{\xi}{|\xi|}=M(\xi)\widehat{f}(\xi).$$
where   $M(\xi)$ is the matrix $(\delta_{k,\ell}-\frac{\xi_k\xi_\ell}{|\xi|^2})_{1\leq k,\ell\leq 3} $.\\
Particularly, if $  u \in \mathcal S(\R^3)^3$, we obtain
$$\ds \mathbb P(  u)_k(x) = \frac{1}{(2\pi)^{\frac 3 2}} \int_{\R^3}  \left( \delta_{kj}-\frac{\xi_k \xi_j}{ \vert \xi \vert^2}\right)
\widehat{  u}_j(\xi) \, e^{i \xi \cdot  x}\,   d\xi,$$
where $\mathcal S(\R^n)$ is the  Schwartz space.
Define also the operator $A_R(D)$ by:
 $$\ds A_R(D)u=\mathbb P J_R(D)u=\mathcal F^{-1}(M(\xi)\chi_{B_R}(\xi)\widehat{u}).$$
%

In which follows, we recall some preliminary results:

\begin{prop}(\cite{HBAF})\label{prop1}\pn
 Let $H$ be a Hilbert space.
\begin{enumerate}
\item The unit ball is weakly compact, that is: if $(x_n)$ is a bounded sequence   in $H$, then there is a subsequence $(x_{\varphi(n)})$ such that
$$(x_{\varphi(n)}|y)\to  (x|y),\;\forall y\in H.$$

\item If $x\in H$ and $(x_n)$   a bounded sequence   in $H$ such that
$\ds\lim_{n\to+\infty}(x_n|y)=  (x|y)$, for all $y\in H,$
then $\|x\|\leq\ds \liminf_{n\to \infty}\|x_n\|.$

\item If $x\in H$ and $(x_n)$ is a bounded sequence   in $H$ such that\\
$\ds\lim_{n\to+\infty}(x_n|y)=  (x|y)$, for all $y\in H$
and
 $\limsup_{n\to \infty}\|x_n\|\leq \|x\|,$
then $\ds \lim_{n\to \infty}\|x_n-x\|=0.$
\end{enumerate}
\end{prop}
We recall   the following product law in the homogeneous Sobolev spaces:

\begin{lem}(\cite{JYC})\label{lem1}\pn
Let $s_1,\ s_2$ be two real numbers and $d\in\N$.

\begin{enumerate}
\item If $s_1<\frac d 2$\; and\; $s_1+s_2>0$, there exists a constant  $C_1=C_1(d,s_1,s_2)$, such that: if $f,g\in \dot{H}^{s_1}(\mathbb{R}^d)\cap \dot{H}^{s_2}(\mathbb{R}^d)$, then $f.g \in \dot{H}^{s_1+s_2-\frac{d}{2}}(\mathbb{R}^d)$ and
$$\|fg\|_{\dot{H}^{s_1+s_2-\frac{d}{2}}}\leq C_1 (\|f\|_{\dot{H}^{s_1}}\|g\|_{\dot{H}^{s_2}}+\|f\|_{\dot{H}^{s_2}}\|g\|_{\dot{H}^{s_1}}).$$

\item If $s_1,s_2<\frac d 2$\; and\; $s_1+s_2>0$ there exists a constant $C_2=C_2(d,s_1,s_2)$ such that: if $f \in \dot{H}^{s_1}(\mathbb{R}^d)$\; and\; $g\in\dot{H}^{s_2}(\mathbb{R}^d)$, then  $f.g \in \dot{H}^{s_1+s_2-\frac{d}{2}}(\mathbb{R}^d)$ and
$$\|fg\|_{\dot{H}^{s_1+s_2-\frac{d}{2}}}\leq C_2 \|f\|_{\dot{H}^{s_1}}\|g\|_{\dot{H}^{s_2}}.$$
\end{enumerate}
 \end{lem}

\begin{lem}\label{lem2}\pn
    Let $ \beta>0$ and $d\in\N$. Then, for all $x,y\in\R^d$, we have
    \begin{equation}\label{eqn-lem2-1}
    \langle |x|^{\beta}x-|y|^{\beta}y ,x-y\rangle\geq \frac{1}{2}(|x|^{\beta}+|y|^{\beta})|x-y|^{2}.
    \end{equation}
\end{lem}
{\bf Proof.} \pn
 Suppose that $|x|>|y|>0$. For $u>v>0$, we have
\begin{equation}\label{eqn-lem2-3}
2\langle ux-vy, x-y\rangle-(u+v)|x-y|^{2}=(u-v)(|x|^2-|y|^2)\geq0.
\end{equation}
It suffices to take $u=|x|^\beta$ and $v=|y|^\beta$, we get the inequality \eqref{eqn-lem2-1}.

 The following result is a generalization of Proposition 3.1 in \cite{J1}.
\begin{prop}\label{prop2}  \pn
Let $\nu_1,\nu_2,\nu_3\in[0,\infty)$, $r_1,r_2,r_3\in(0,\infty)$ and $f^0\in L^2_\sigma(\R^3)$. \\
For $n\in\N$, let $F_n:\R^+\times\R^3\to \R^3$ be a measurable function in $C^1(\R^+,L^2(\R^3))$ such that $$A_n(D)F_n=F_n,\;F_n(0,x) =A_n(D)f^0(x)$$ and

\begin{enumerate}
\item [(E1)]
$\ds  \partial_t F_n+\sum_{k=1}^3\nu_k|D_k|^{2r_k} F_n+ A_n(D){\rm div}\,(F_n\otimes F_n)+ A_n(D)h(|F_n|)F_n =0.$

\item [(E2)]
\beqq
&&\ds \|F_n(t)\|_{L^2}^2+2\sum_{k=1}^3\nu_k\int_0^t\||D_k|^{r_k} F_n(s)\|_{L^2}^2ds\\
&&\hskip 2cm +2 a\int_0^t\|h(|F_n(s)|)|F_n(s)|^2\|_{L^1}ds \leq \|f^0\|_{L^2}^2.
\eeqq
\end{enumerate}
 where  $\ds h(z)=\alpha z^{\beta-1},$ with  $\alpha>0$ and $\beta>3$.
  Then: for every $\varepsilon>0$ there is $\delta=\delta(\varepsilon,\alpha,\beta,\nu_1,\nu_2,\nu_3,r_1,r_2,r_3,\|f^0\|_{L^2})>0$
  such that: for all $t_1,t_2\in\R^+$, we have

\begin{equation}\label{eqn-1}
\Big(|t_2-t_1|<\delta\Longrightarrow \|F_n(t_2)-F_n(t_1)\|_{H^{-s_0}}<\varepsilon\Big),\;\forall n\in\N,
\end{equation}
with $\ds  s_0\ge \max(3,2r_1,2r_2,2r_3).$
\end{prop}

{\bf Proof.} The proof is similar to that of Proposition 2.4 in \cite{MJ1}.

\section{\bf Proof of Theorem \ref{th2}}

\subsection{Existence of Solution}\ \\
 Consider the approximate system: 
$$(NSD_{n})
  \begin{cases}
     \partial_t u
 -\Delta J_nu+ J_n(J_nu.\nabla J_nu)  + \alpha J_n[ |J_nu|^{\beta-1}J_nu] =\;\;-\nabla p_n\hbox{ in } \mathbb R^+\times \mathbb R^3\\
 p_n=(-\Delta)^{-1}\Big({\rm div}\,J_n(J_nu.\nabla J_nu)  + \alpha {\rm div}\,J_n[|J_nu|^{\beta-1}J_nu]\Big)\\
     {\rm div}\, u = 0 \hbox{ in } \mathbb R^+\times \mathbb R^3\\
    u(0,x) =J_nu^0(x) \;\;\hbox{ in }\mathbb R^3,
  \end{cases}
$$
where $J_n$ is the Friedritch operator defined by
$\ds J_n(D)f=\F^{-1}(\chi_{B_n} \widehat{f}) $ and
 $B_n$   the ball of center $0$ and radius $n\in\N$.
\begin{enumerate}
\item[$\bullet$] By Cauchy-Lipschitz theorem, there exists a unique solution $u_n\!\in\! C^1(\R^+,L^2_\sigma(\R^3))$ of the system $(NSD_{n})$ such that $J_nu_n\!=\!u_n$ and

\begin{equation}\label{eq-energyn}\|u_n(t)\|_{L^2}^2+2\int_0^t\|\nabla u_n(s)\|_{L^2}^2ds +2\alpha\int_0^t\|u_n(s)\|_{L^{\beta+1}}^{\beta+1}ds\leq \|u^0\|_{L^2}^2.\end{equation}

\item[$\bullet$]
The sequence $(u_n)_n$ is bounded in $L^2(\R^3))$ and on $H^{1}(\R^3)$.\\
 Using proposition \ref{prop2} and the  interpolation method, we deduce that the sequence $(u_n)_n$ is equicontinuous on $H^{-1}(\R^3)$.

\item[$\bullet$] 
For  $(T_q)_q $   a strictly  increasing sequence such that  $\ds\lim_{q\to+\infty} T_q=\infty$, consider a sequence
of functions $(\theta_q)_{q }$ in $C_0^\infty(\R^3)$  such that

$$\left\{\begin{array}{l}
\theta_q(x)=1,\ {\rm for}\  |x|\le   q+\frac{5}{4}\\
\theta_q(x)=0,\ {\rm for}\   |x|\ge  q+2 \\
0\leq \theta_q\leq 1.
\end{array}\right.$$

Using the energy estimate  \eqref{eq-energyn}, the equicontinuity of the sequence $(u_n)_n$   on $H^{-1}(\R^3)$
 and classical argument by combining Ascoli's theorem and the Cantor diagonal process, there exists a subsequence $(u_{\varphi(n)})_n$   and
$u\in L^\infty(\R^+,L^2(\R^3))\cap C(\R^+,H^{-3}(\R^3))$ such that: for all $q\in\N$,

\begin{equation}\label{eq-cv}\lim_{n\to \infty}\|\theta_q(u_{\varphi(n)}(t)-u(t))\|_{L^\infty([0,T_q],H^{-4})}=0.\end{equation}
In particular, the sequence $(u_{\varphi(n)}(t))_n$ converges weakly in $L^2(\R^3)$ to $u(t)$ for all
$t\geq0$.

\item[$\bullet$] Using the same method in \cite{J1}, we obtain:
\begin{equation}\label{eqn-23}
\|u(t)\|_{L^2}^2\!+\!2\int_0^t\!\|\nabla u(s)\|_{L^2}^2ds\!+\!2\alpha\int_0^t\!\|u(s)\|_{L^{\beta+1}}^{\beta+1}ds\!\leq\! \|u^0\|_{L^2}^2.
\end{equation}
for all $t\geq0$, and $u$ is a solution of the system $(NSD)$.
\end{enumerate}

\subsection{Continuity of the solution in $L^2$}\pn
   By the inequality (\ref{eqn-23}), we have
 $\ds \limsup_{t\to 0}\|u(t)\|_{L^2}\leq\|u^0\|_{L^2} $ and using
proposition \ref{prop1}-(3), we get
 $\ds \limsup_{t\to 0}\|u(t)-u^0\|_{L^2}=0.$
This ensures the continuity of the solution  $u$ at $0$. To prove the continuity on $\R$, consider the  functions
 $\ds v_{n,\varepsilon}(t)=u_{\varphi(n)}(t+\varepsilon),\;p_{n,\varepsilon}(t)=p_{\varphi(n)}(t+\varepsilon),$
for $n\in\N$ and $\varepsilon>0$. We have:

\beqq
\partial_tu_{\varphi(n)}-\Delta u_{\varphi(n)}+J_{\varphi(n)}(u_{\varphi(n)}.\nabla u_{\varphi(n)})
+\alpha J_{\varphi(n)}(|u_{\varphi(n)}|^{\beta-1} u_{\varphi(n)})&=&-\nabla p_{\varphi(n)} \\
\partial_tv_{n,\varepsilon}-\Delta v_{n,\varepsilon}+J_{\varphi(n)}(v_{n,\varepsilon}.\nabla
v_{n,\varepsilon})+\alpha J_{\varphi(n)}(|v_{n,\varepsilon}|^{\beta-1} v_{n,\varepsilon})&=&-\nabla p_{n,\varepsilon}.
\eeqq
The function $w_{n,\varepsilon}=u_{\varphi(n)}-v_{n,\varepsilon}$ fulfills the following:

\beqq
&&\partial_tw_{n,\varepsilon}-\Delta w_{n,\varepsilon}  +\alpha J_{\varphi(n)}\Big(|u_{\varphi(n)}|^{\beta-1}
 u_{\varphi(n)}-|v_{\varphi(n)}|^{\beta-1}  v_{n,\varepsilon}\Big)\\
&&\hskip  3cm = -\nabla (p_{\varphi(n)}-p_{n,\varepsilon})+J_{\varphi(n)}(w_{n,\varepsilon}.\nabla w_{n,\varepsilon})\\
&&\hskip  3cm-J_{\varphi(n)}(w_{n,\varepsilon}.\nabla u_{\varphi(n)})
 - J_{\varphi(n)}(u_{\varphi(n)}.\nabla w_{n,\varepsilon}).
 \eeqq
Taking the scalar product with $w_{n,\varepsilon}$ in $L^2(\R^3)$  and using the fact that
$\langle w_{n,\varepsilon}.\nabla w_{n,\varepsilon},w_{n,\varepsilon}\rangle=0$ and ${\rm div}\ w_{n,\varepsilon}=0$, we get

\begin{eqnarray}\label{eqn-24}
\frac{1}{2}\frac{d}{dt}\|w_{n,\varepsilon}(t)\|_{L^2}^2+\|\nabla
w_{n,\varepsilon}(t)\|_{L^2}^2&
 +&\alpha \langle J_{\varphi(n)}\Big( |u_{\varphi(n)}|^{\beta-1}u_{\varphi(n)}- |v_{n,\varepsilon}|^{\beta-1}v_{n,\varepsilon}\Big);w_{n,\varepsilon}\rangle_{L^2}
\nonumber\\
&=& -\langle J_{\varphi(n)}(w_{n,\varepsilon}.\nabla u_{\varphi(n)});w_{n,\varepsilon}\rangle _{L^2} .
\end{eqnarray}
By inequality  \eqref{eqn-lem2-1}, we have

\begin{eqnarray}
&&\langle J_{\varphi(n)}\Big( |u_{\varphi(n)}|^{\beta-1}u_{\varphi(n)}-( |v_{n,\varepsilon}|^{\beta-1}v_{n,\varepsilon}\Big);
w_{n,\varepsilon}\rangle _{L^2}\nonumber\\
 &&\hskip 5cm=\langle ( |u_{\varphi(n)}|^{\beta-1} u_{\varphi(n)}-
 |v_{n,\varepsilon}|^{\beta-1} v_{n,\varepsilon};J_{\varphi(n)}w_{n,\varepsilon}\rangle _{L^2}\nonumber\\
 &&\hskip 5cm= \langle ( |u_{\varphi(n)}|^{\beta-1}u_{\varphi(n)}- |v_{n,\varepsilon}|^{\beta-1}v_{n,\varepsilon};w_{n,
\varepsilon}\rangle_{L^2}\nonumber\\
 && \hskip 5cm\geq  \frac{1}{2}\int_{\R^3}\Big( |u_{\varphi(n)}|^{\beta-1}+ |v_{n,\varepsilon}|^{\beta-1}\Big)|w_{n,\varepsilon}|^2\nonumber\\
&&\hskip 5cm  \geq    \frac{1}{2}\int_{\R^3}
|u_{\varphi(n)}|^{\beta-1}|w_{n,\varepsilon}|^2,\nonumber
\end{eqnarray}
which implies
{\footnotesize
\begin{equation}\label{eqn41}\alpha\langle J_{\varphi(n)}\Big( |u_{\varphi(n)}|^{\beta-1}u_{\varphi(n)}-( |v_{n,\varepsilon}|^{\beta-1}v_{n,\varepsilon}\Big);
w_{n,\varepsilon}\rangle _{L^2}\geq \frac{\alpha}{2}\int_{\R^3}
|u_{\varphi(n)}|^{\beta-1}|w_{n,\varepsilon}|^2.\end{equation}}
Also, we have
\begin{eqnarray}\label{eqn42}
|\langle J_{\varphi(n)}(w_{n,\varepsilon}.\nabla u_{\varphi(n)});w_{n,\varepsilon}\rangle _{L^2}|
&\leq&\ds \int_{\R^3}|w_{n,\varepsilon}|.|u_{\varphi(n)}|.|\nabla w_{n,\varepsilon}| \nonumber\\
&\leq&\ds \frac{1}{2}\int_{\R^3}|w_{n,\varepsilon}|^2|u_{\varphi(n)}|^2 +\frac{1}{2}\|\nabla w_{n,\varepsilon}\|_{L^2}^2\nonumber.
\end{eqnarray}
By using the convex inequality
$$ab\leq \frac{a^p}{p}+\frac{b^q}{q}\leq a^p+b^q$$
with $p = \ds\frac{\beta-1}{2},\ \ q=\ds\frac{\beta-1}{\beta-3},\ \
a=\ds|w_{n,\varepsilon}|^2(\frac \alpha 2)^{\frac{2}{\beta-1}},\ \
b=\ds(\frac 2\alpha)^{\frac{2}{\beta-1}},$

we get
\begin{eqnarray}\label{eqn42}
|\langle J_{\varphi(n)}(w_{n,\varepsilon}.\nabla u_{\varphi(n)});w_{n,\varepsilon}\rangle _{L^2}|
&\leq&\ds \frac{\alpha}{4}\int_{\R^3}|w_{n,\varepsilon}|^{\beta-1}|u_{\varphi(n)}|^2+C_{\alpha,\beta}\|w_{n,\varepsilon}\|_{L^2}^2+\frac{1}{2}\|\nabla w_{n,\varepsilon}\|_{L^2}^2,\nonumber
\end{eqnarray}
with $C_{\alpha,\beta}=\frac{1}{2}(\frac 2\alpha)^{\frac{2}{\beta-3}}$. Combining this inequality and inequalities \eqref{eqn-24}, \eqref{eqn41} and  \eqref{eqn42}, we get

$$ \frac{d}{dt}\|w_{n,\varepsilon}\|_{L^2}^2+ \|\nabla w_{n,\varepsilon}\|_{L^2}^2\leq   2C_{\alpha,\beta} \|w_{n,\varepsilon}\|_{L^2}^2.$$
By Gronwall Lemma, we deduce the following:
$$\|w_{n,\varepsilon}(t)\|_{L^2} \leq  \|w_{n,\varepsilon}(0)\|_{L^2} e^{C_{\alpha,\beta}t},$$
and
$$\|u_{\varphi(n)}(t+\varepsilon)-u_{\varphi(n)}(t)\|_{L^2} \leq
 \|u_{\varphi(n)}(\varepsilon)-u_{\varphi(n)}(0)\|_{L^2} e^{C_{\alpha,\beta}t}.$$
For $t_0>0$ and $\varepsilon\in(0,t_0)$, we have

$$\|u_{\varphi(n)}(t_0+\varepsilon)-u_{\varphi(n)}(t_0)\|_{L^{2}}\leq\|u_{\varphi(n)}(\varepsilon)-u_{\varphi(n)}(0)\|_{L^{2}} e^{C_{\alpha,\beta}t_0}.$$

$$\|u_{\varphi(n)}(t_0-\varepsilon)-u_{\varphi(n)}(t_0)\|_{L^{2}}\leq\|u_{\varphi(n)}(\varepsilon)-u_{\varphi(n)}(0)\|_{L^{2}} e^{C_{\alpha,\beta}t_0}.$$
So
\beqq
\| u_{\varphi(n)}(\varepsilon)-u_{\varphi(n)}(0) \|_{L^2}^2 &=&
\| J_{\varphi(n)} u_{\varphi(n)}(\varepsilon)-J_{\varphi(n)}u_{\varphi(n)}(0) \|_{L^2}^2 \\
&=&\|J_{\varphi(n)}\left(u_{\varphi(n)}(\varepsilon)  -u^0\right)\|_{L^2}^2\\
&\le&\| u_{\varphi(n)}(\varepsilon)-u^0 \|_{L^2}^2\\
&\le& \| u_{\varphi(n)}(\varepsilon)\|_{L^2}^2+\|u^0 \|_{L^2}^2-2Re\langle u_{\varphi(n)}(\varepsilon),u^0\rangle\\
&\le& 2\|u^0 \|_{L^2}^2-2Re\langle u_{\varphi(n)}(\varepsilon),u^0\rangle.
\eeqq
But $\ds \lim_{n\to+\infty}\langle u_{\varphi(n)}(\varepsilon),u^0\rangle=\langle u (\varepsilon),u^0\rangle$, hence

$$\liminf_{n\to \infty}\|u_{\varphi(n)}(\varepsilon)-u_{\varphi(n)}(0)\|^{2}_{L^{2}}
\leq 2\|u^0\|^{2}_{L^{2}}-2Re\langle
u(\varepsilon);u^0\rangle_{L^2}.$$ Moreover, for all $q,N\in\N$
\beqq
\|J_N\Big(\theta_q.(u_{\varphi(n)}(t_0\pm\varepsilon)-u_{\varphi(n)}(t_0))\Big)\|^{2}_{L^2}
 &\leq & \|\theta_q.(u_{\varphi(n)}(t_0\pm\varepsilon)-u_{\varphi(n)}(t_0))\|^{2}_{L^2}\\
 &\leq&
 \|u_{\varphi(n)}(t_0\pm\varepsilon)-u_{\varphi(n)}(t_0)\|^{2}_{L^2}.
 \eeqq

 Using  (\ref{eq-cv}) we get, for $q$ big enough,
 $$\|J_N\Big(\theta_q.(u(t_0\pm\varepsilon)-u(t_0))\Big)\|_{L^2}
 \leq \liminf_{n\to \infty}\|u_{\varphi(n)}(t_0\pm\varepsilon)-u_{\varphi(n)}(t_0)\|_{L^2}.$$
 Then
$$\|J_N\Big(\theta_q.(u(t_0\pm\varepsilon)-u(t_0))\Big)\|^{2}_{L^2}
 \leq 2\Big(\|u^0\|^{2}_{L^{2}}-Re\langle u(\varepsilon);u^0\rangle_{L^2}\Big)e^{2C_{\alpha,\beta}t_0}.$$
By applying the Monotone Convergence Theorem in the order $N $ and $q $, we get
$$\|u(t_0\pm\varepsilon)-u(t_0))\|^{2}_{L^2}
 \leq 2\Big(\|u^0\|^{2}_{L^{2}}-Re\langle u(\varepsilon);u^0\rangle_{L^2}\Big)e^{ 2C_{\alpha,\beta}t_0}.$$
Using the continuity at 0 and make $\varepsilon\to 0$, we get the continuity at $t_0$.

\subsection{Uniqueness}\ \\
Let $u,v$ be two solutions of $(NSD)$ in the space

$$C_b(\R^+,L^2(\R^3))\cap L^2(\R^+,\dot H^1(\R^3))\cap L^{\beta+1}(\R^+,L^{\beta+1}(\R^3)).$$
The function $w=u-v$ satisfies the following:

$$\partial_tw-\Delta w+\alpha \Big( |u|^{\beta-1}u- |v|^{\beta-1} v\Big)= -\nabla (p-\tilde p)+w.\nabla w-w.\nabla u-
u.\nabla w.$$
Taking the scalar product in $L^2$ with $w$, we get

$$\frac{1}{2}\frac{d}{dt}\|w\|_{L^2}^2+\|\nabla w\|_{L^2}^2+\alpha
\langle \Big( |u|^{\beta-1} u- |v|^{\beta-1}v\Big);w\rangle _{L^2}=-\langle w.\nabla u;w\rangle _{L^2}.$$
By adapting the same method for the proof of the continuity of such solution in $L^2(\R^3)$, with
$u, v,w $ instead of $u_{\varphi(n)}, v_{n,\varepsilon}, w_{n,\varepsilon}$ in oreder,

we find

$$
\alpha\langle \Big( |u|^{\beta-1}u- |v|^{\beta-1} v\Big);w\rangle _{L^2}\geq
\frac{\alpha}{2}\int_{\R^3}|u|^{\beta-1}|w|^2.
$$
and
$$
|\langle w.\nabla u;w\rangle _{L^2}|\leq \frac{\alpha}{4}\int_{\R^3}|w|^2|u|^2+C_{\alpha,\beta}\|w\|_{L^2}^2+\frac{1}{2}\|\nabla w\|_{L^2}^2.$$
Combining the above inequalities, we find the following energy estimate:
$$ \frac{d}{dt}\|w(t)\|_{L^2}^2+ \|\nabla w(t)\|_{L^2}^2 \leq  2C_{\alpha,\beta}\|w(t)\|_{L^2}^2.$$
By Gronwall Lemma, we obtain
$$\|w(t)\|_{L^2}^2+\int_0^t\|\nabla w\|_{L^2}^2 \leq \|w^0\|_{L^2}^2e^{2C_{\alpha,\beta}t}.$$
As $w^0=0$, then $w=0$ and $u=v$, which implies the uniqueness.

\subsection{Asymptotic Study of the Global Solution}\pn
To prove the asymptotic behavior  \eqref{eqth2-2},   we need some preliminaries lemmas:
\begin{lem}\label{lem1}\pn
 If $u$ is a global solution of {$(NSD)$} with $\beta\!\geq\!\frac{10}3$, then  $u \in  L^{\beta}(\R^+\times\R^3)$.
\end{lem}
\noindent{\bf Proof.}\pn
Let  $E_1=\{(t,x):\ |u(t,x)|\leq 1\}$ and $E_2=\{(t,x):\ |u(t,x)|> 1\},$
$\ds  L_1=\ds\int_{E_1}|u(s,x)|^{\beta}dxds $  and $\ds L_2=\ds\int_{E_2}|u(s,x)|^{\beta}dxds.$
 We  have
%
%

 \beqq
L_1&=&\ds\int_{E_1}|u(s,x)|^{\beta}dxds =\ds\int_{E_1}|u(s,x)|^{\beta-\frac{10}3}|u(s,x)|^{\frac{10}3}dxds \\
& \leq& \ds \int_{0}^\infty\|u(s)\|_{L^{\frac{10}3}}^{\frac{10}3}ds .
\eeqq
By using the Sobolev injection
$\dot H^{\frac{3}5}(\R^3)\hookrightarrow L^{\frac{10}3}(\R^3)$, we get
\begin{equation}\label{eqasym1}
L_1\leq C \int_{0}^\infty\|u(s)\|_{\dot H^{\frac{3}5}}^{\frac{10}3}ds.
\end{equation}

By interpolation inequality
$\ds \|u(s)\|_{\dot H^{\frac{3}5}}\leq \|u(s)\|_{\dot H^0}^{\frac{2}5}\|u(s)\|_{\dot H^1}^{\frac{3}5}$,
we obtain
\begin{equation}\label{eqasym2}
L_1\leq C \int_{0}^\infty\|u(s)\|_{L^2}^{\frac{4}3}\|\nabla u(s)\|_{L^2}^2\leq C
\|u^0\|_{L^2}^{\frac{4}3}\int_{0}^\infty\|\nabla u(s)\|_{L^2}^2.\end{equation}
For the therm $L_2$, we have

$$L_2=\ds\int_{X_2}|u(s,x)|^\beta dxds\le \ds\int_{0}^\infty\int_{\R^3}|u(s,x)|^{\beta+1}dxds.
$$

Hence
$$\|u\|_{L^{\beta}(\R^+\times\R^3)}\leq C
\|u^0\|_{L^2}^{\frac{4}3}\int_{0}^\infty\|\nabla
u(s)\|_{L^2}^2ds +\int_{0}^\infty\int_{\R^3}|u(s,x)|^{\beta+1}dxds.$$
Therefore $u\in L^{\beta}(\R^+\times\R^3)$.

\begin{lem}\pn
 If $u$ is a global solution of {$(NSD)$}, with $\beta\!\geq\! \frac{10}3$ , then  $\ds \lim_{t\!\to\!\infty}\|u(t)\|_{H^{\!-\!2}}\!=\!0$.
\end{lem}
\noindent{\bf Proof.}\pn
For $\varepsilon>0$, using  the energy inequality
(\ref{eqth2-1}) and Lemma \ref{lem1}, there exists   $t_0\geq0$ such that

\begin{equation}\label{asym.eq1}
\|\nabla u\|_{L^2([t_0,\infty)\times\R^3)}<\frac{\varepsilon}{4},
\end{equation}
and
\begin{equation}\label{asym.eq2}
\|u\|_{L^\beta([t_0,\infty)\times\R^3)}<\frac{\varepsilon}{4}.
\end{equation}
Now, consider the following system
\begin{equation}\label{$4.6$}
 \left\{ \begin{matrix}
     \partial_t v
 -\nu\Delta v+ v.\nabla v  +\alpha|v|^{\beta-1}v =\;\;-\nabla q \hfill&\hbox{ in } \mathbb R^+\times \mathbb R^3\\
     {\rm div}\, v= 0 \hfill&\hbox{ in } \mathbb R^+\times \mathbb R^3\\
    v(0,x) =u(t_0,x) \;\;\hfill&\hbox{ in }\mathbb R^3\hfill&.
\end{matrix}\right. \tag{$NSD'$}
\end{equation}
By the existence and uniqueness part, the system ($NSD'$) has a
unique global solution $v\in C_b(\R^+,L^2(\R^3))\cap L^2(\R^+,\dot H^1(\R^3))\cap L^{\beta+1}(\R^+,L^{\beta+1}(\R^3))$
such that $v(t_0)=u(t_0,x)$ and $q(t)=p(t_0+t).$
The energy estimate  for this system is as follows:
$$\|v(t)\|_{L^2}^2+2\int_0^t\|\nabla
v(s)\|_{L^2}^2ds +2a\int_0^t\|v(s)\|_{L^{\beta+1}}^{\beta+1}\leq\|u(t_0)\|_{L^2}^2\leq
\|u^0\|_{L^2}^2.$$
By  the  Duhamel formula, we obtain
$$\ds v(t,x)=e^{t\Delta}v^0(x)+f(t,x)+g(t,x),$$ where
$$f(t,x)=-\int_0^te^{(t-s)\Delta}\mathbb P{\rm div\,}(v\otimes v)(s,x)ds$$
and
$$g(t,x)=-\alpha\int_0^te^{(t-s)\Delta}\mathbb P{\rm div\,}|v(s,x)|^{\beta-1}v(s,x)ds.$$
By   Dominated Convergence Theorem, 
 $\ds \lim_{t\rightarrow\infty}\|e^{t\Delta}v^0\|_{L^2}=0$ and hence $\ds  \lim_{t\rightarrow\infty}\|e^{t\Delta}v^0\|_{H^{-2}}=0.$\\
 Moreover,
\beqq
\|f(t)\|_{H^{-2}}^2&\leq&\ds\|f(t)\|_{H^{-\frac{1}2}}^2 \leq\ds\|f(t)\|_{\dot H^{-\frac{1}2}}^2\\
&\leq&\ds\int_{\R^3}|\xi|^{-1}\left(\int_0^te^{-(t-s)|\xi|^2}|\mathcal F{\rm div}(v\otimes v)(s,\xi)|ds\right)^2d\xi\\
&\leq&\ds\int_{\R^3}|\xi|\left(\int_0^te^{-(t-s)|\xi|^2}|\mathcal F(v\otimes v)(s,\xi)|ds\right)^2d\xi.
\eeqq
Since
\beqq
\ds\left(\int_0^te^{-(t-s)|\xi|^2}|\mathcal F(v\otimes v)(s,\xi)|ds\right)^2
&\leq&\ds\left(\int_0^te^{-2(t-s)|\xi|^2}ds\right) \int_0^t|\mathcal F(v\otimes v)(s,\xi)|^2ds \\
&\leq&\ds|\xi|^{-2} \int_0^t|\mathcal F(v\otimes v)(s,\xi)|^2ds,
\eeqq
then 
\beqq
\|f(t)\|_{H^{-2}}^2dt&\leq&\ds\int_{\R^3}|\xi|^{-1}\int_0^t|\mathcal F(v\otimes v)(s,\xi)|^2dsd\xi \\
&\leq&\ds\int_0^t(\int_{\R^3}|\xi|^{-1}|  (v\otimes v)(s,\xi)|^2d\xi)ds=\ds\int_0^t\|v\otimes v)(s)\|_{\dot H^{-\frac{1}2}}^2ds.
\eeqq

\noindent
Using the  product law in homogeneous Sobolev spaces, with
$s_1=0$, $s_2=1$, we get
\beqq
\|f(t)\|_{H^{-2}}^2dt&\leq&\ds C\int_0^t\|v(s)\|_{L^2}^2\|\nabla
v(s)\|_{L^2}^2ds.
\eeqq
Using inequalities  \eqref{asym.eq1}  and \eqref{asym.eq2}, we get
\beqq
\|f(t)\|_{H^{-2}}^2dt&\leq&\ds C\|u^0\|_{L^2}^2\int_0^t\|\nabla
u(t_0+s)\|_{L^2}^2ds\\
&\leq&\ds C\|u^0\|_{L^2}^2\int_0^\infty\|\nabla u(t_0+s)\|_{L^2}^2ds\\
&\leq&\ds C\|u^0\|_{L^2}^2\int_{t_0}^\infty\|\nabla u(s)\|_{L^2}^2ds\\
&\leq&\ds C\|u^0\|_{L^2}^2\frac{\varepsilon^2}{9(C\|u^0\|_{L^2}^2+1)},
\eeqq
which implies that
$$\|f(t)\|_{H^{-2}}<\frac \varepsilon3,\;\forall t\geq 0.$$

For an estimation of $\|g(t)\|_{H^{-2}}$ and using
$$L^1(\R^3)\hookrightarrow H^{-s}(\R^3),\;\forall s>3/2,$$
with $s=2$, we get
\beqq
\|g(t)\|_{H^{-2}}^2dt&\leq&\ds
\int_{\R^3}(1+|\xi|^2)^{-2}\left(\int_0^te^{-(t-s)|\xi|^2}|\mathcal
F(|v |^{\beta-1}v)(s,\xi)|ds\right)^2d\xi\\
&\leq&\ds
C\left(\int_0^t\|(|v |^{\beta-1}v)(s,.)\|_{L^1(\R^3)}ds\right)^2\\
&\leq&\ds
C\left(\int_0^t\||v(s,.)|^{\beta}\|_{L^1(\R^3)}ds\right)^2\\
&\leq&\ds C\|v\|_{L^{\beta}(\R^+\times\R^3)}^2,
\eeqq
where $\ds C=\int_{\R^3}(1+|\xi|^2)^{-2}d\xi.$\\
Also  using inequality  \eqref{asym.eq2}, we get
\beqq
\|g(t)\|_{H^{-2}}^2dt&\leq&C\||u(t_0+.)\|_{L^{\beta}(\R^+\times\R^3)}^2\\
&\leq&C\|u\|_{L^\beta([t_0,\infty)\times\R^3)}^2
 \leq C\frac{\varepsilon^2}{9C}.
\eeqq

which implies that
$\ds \|g(t)\|_{H^{-2}}<\frac \varepsilon 3,\;\forall t\geq 0.$

Combining the above inequalities, we obtain

$$\lim_{t\rightarrow\infty}\|u(t)\|_{H^{-2}}=0.$$

\begin{lem}\pn
 If $u$ is a global solution of $(NSD)$ and $\beta\geq\frac{10}3$, then  $\ds \lim_{t\rightarrow\infty}\|u(t)\|_{L^2}=0.$
\end{lem}
\noindent{\bf Proof.}\pn
  Let  
$$w_1 = {\bf 1}_{|D|<1}u=\mathcal F^{-1}\big({\bf 1}_{|\xi|<1}\widehat{u}\big)\quad{\rm and}\quad
w_2 = {\bf 1}_{|D|\geq1}u=\mathcal F^{-1}\big({\bf
1}_{|\xi|\geq1}\widehat{u}\big).$$
Using   the second step, we get
$$\|w_1(t)\|_{L^2}=c_0\|w_1(t)\|_{H^0}\leq 2c_0\||w_1(t)\|_{H^{-2}}\leq 2\||u(t)\|_{H^{-2}},$$
which implies
$$\lim_{t\rightarrow\infty}\|w_1(t)\|_{L^2}=0.$$
 For $\varepsilon>0$, there is a   $t_1>0$ such that
$$\|w_1(t)\|_{L^2}<\frac \varepsilon2,\;\forall t\geq t_1.$$
We have
$$\int_{t_1}^\infty\|w_2(t)\|_{L^2}^2dt\leq \int_{t_1}^\infty\|\nabla w_2(t)\|_{L^2}^2dt\leq \int_{t_1}^\infty\|\nabla u(t)\|_{L^2}^2dt<\infty.$$
Since the map $t\longmapsto \|w_2(t)\|_{L^2}$   is continuous, there exists  $t_2\geq t_1$ such that
 $\ds \|w_2(t_2)\|_{L^2}<\frac \varepsilon2.$  Hence  
$$\|u(t_2)\|_{L^2}^2=\|w_1(t_2)\|_{L^2}^2+\|w_2(t_2)\|_{L^2}^2<\frac {\varepsilon^2}2.$$
Using the following energy estimate
$$
\|u(t)\|_{L^2}^2+2\int_{t_2}^t\|\nabla
u(s)\|_{L^2}^2ds+2\alpha\int_{t_2}^t\|u(s)\|_{L^{\beta+1}}ds\leq\|u(t_2)\|_{L^2}^2,\,\forall t\geq t_2,$$ we get
$$\|u(t)\|_{L^2}<\varepsilon,\;\forall t\geq t_2,$$
and the proof is completed.

%

\end{document}